\documentclass[%
aip,
sd,%
amsmath,amssymb,
reprint,%
author-numerical,%
]{revtex4-1}

\draft 

\usepackage[T1]{fontenc}
\usepackage[utf8]{inputenc}
\usepackage{amsmath}
\usepackage[american]{babel}

\usepackage{booktabs}
\usepackage[small,bf]{caption}
\usepackage{graphicx}
\usepackage{csquotes}
\usepackage{subfigure}

\usepackage[normalem]{ulem}

\usepackage{xparse}

\NewDocumentCommand{\evalat}{sO{\big}mm}{
	\IfBooleanTF{#1}
	{\mleft. #3 \mright|_{#4}}
	{#3#2|_{#4}}%
}

\newcommand{\dd}[1]{\mathrm{d}#1}

\newcommand\Rey{\mbox{\textit{Re}}}

\newcommand{\pder}[2]{\frac{\partial#1}{\partial#2}} 
\newcommand{\tder}[2]{\frac{\dd#1}{\dd#2}} 

\begin{document}


\title{
Analytic investigation of the compatibility condition and the initial evolution of a smooth velocity field for the Navier-Stokes equation in a channel configuration}   



\author{P\'eter Tam\'as Nagy}
\email{pnagy@hds.bme.hu}
\author{Gy\"{o}rgy Pa\'{a}l}%
\email{paal@hds.bme.hu}
\affiliation{Department of Hydrodynamic Systems, Faculty of Mechanical Engineering,Budapest University of Technology and Economics, Budapest, H-1111, Hungary}


\date{\today}

\begin{abstract}
	A partial differential equation has usually a regular solution at the initial time if the initial condition is smooth in space, fulfills the governing equations and is compatible with the boundary condition. In the case of Navier-Stokes equation, the initial velocity field must also be divergence--free. It is common belief that the initial condition is compatible with the boundary condition if the initial condition fulfills the boundary condition but this is not sufficient.  Such a field does not necessarily fulfill the full compatibility condition of the Navier-Stokes equation. If the condition is violated, the solution is not regular at the initial time. This issue has been known for a while but not in the full breadth of the fluid dynamics community. In this paper, a practical calculation method is presented for checking the compatibility condition. Furthermore, a smooth initial condition is presented in a periodic channel flow that violates the condition and has no spatially smooth solution at initial time. The calculations were were performed in an analytical framework.
	
	The results for a channel configuration show that in the absence of wall--normal velocity the condition is always fulfilled and the problem has a regular solution. If the wall--normal velocity component is non--zero, the condition is usually not fulfilled but with optimization methods counter-examples can be generated. The generation procedure of such a field is useful to provide correct initial conditions for sensitive time-dependent numerical simulations.  
	 The presented methods can provide insight about the applicability of the chosen initial conditions.
\end{abstract}

\maketitle 

\section{Introduction}
Some mathematical proofs about the existence of a smooth Navier-Stokes solution were derived in the last century. The first theorems were summarized in the book of \citet{Ladyzhenskaia1969}. The Navier-Stokes differential equations have a unique solution in two dimensions until arbitrary time. Unfortunately, this statement has not been proved yet in three dimensions 
but it is known that the unique solution must exist for a certain time $T^*$ on unbounded or periodic domains. However, if the domain is bounded, a wall is present, the initial condition cannot be arbitrary. This problem was briefly discussed in the aforementioned book, and it was later investigated thoroughly by \citet{Temam1982}. He proved that the solution is regular in three-dimensional space at $t=0$, if and only if the initial condition fulfills the compatibility condition. 
However, this condition is not well-known in the fluid dynamics community and therefore will be presented in Section \ref{sec_method}. Later, he gave a physical explanation for these mathematical results \citep{Temam2006}. The problem is that if the velocity is prescribed at the boundaries (Dirichlet boundary condition), the calculation of pressure leads to an overdetermined problem. The usual method is using the wall-normal momentum equation to obtain a Neuman boundary condition for the pressure. Still, the calculated pressure can violate the tangential component of the momentum equation in an arbitrary case. Using the tangential component of the momentum equation to obtain a boundary condition would lead to different pressure field \citep{Orszag1974}, unless the compatibility condition is fulfilled.  
The same problem appeared in computational fluid mechanics (CFD), as discussed by \citet{Gresho1987}. The main question asked by the authors was, how the boundary condition for the pressure should be handled. If the compatibility condition is not fulfilled, the Dirichlet and the Neumann boundary condition lead to different solutions. They recommend the usage of Neumann boundary condition at the initial time. Furthermore, they propose multiple methods to handle the problems if the condition is not fulfilled, which is generally the case in CFD simulations. In addition, they showed numerical examples. A similar technique was adapted to solve efficiently the incompressible Navier-Stokes equation using the pressure Poisson equation with a Neumann boundary condition by \citet{Johnston2004}.
\citet{Gallavotti2002} discussed the problem with the initial condition in Chapter 2.1 of his book. He suggests a theoretically possible numerical procedure to handle the issue by extending the domain with a thin additional volume. According to his solution, the velocity is not zero at the original boundary but it is reduced by an extra "friction" term in the thin layer over the boundary. 
Nevertheless, he admits that the numerical solution of the problem using the suggested method is challenging. 
It must be mentioned that these techniques were developed to obtain numerical results, but according to the theorem of \citet{Temam1982}, the regular solution does not exist in these cases at $t=0$.
This problem can cause errors at the beginning of simulation, especially if the initial behavior of the flow is critical. Another example is the energy stability analyis of fluid flows, namely the Reynolds-Orr\citep{Orr1907} equation. In this case, the stability problem is rewritten to a variational problem. The perturbation velocity field is varied to maximize the temporal growth rate of the kinetic energy. The base flow is stable if this maximum is below zero. However, the velocity field usually does not fulfill the compatibility condition and the using it as an initial condition may lead to nonphysical results, at least in the beginning of the simulation. 

Although the Navier-Stokes equation rarely has an analytical solution, there are many counterexamples. They are derived usually with some approximations or assumptions. 
In addition, some mathematical papers have discussed the extraordinary analytical solutions of the Navier-Stokes equation in the weak form. For example, \citet{Heywood1980} showed that multiple solutions exist in the case of flow through a hole. He pointed out the need for a further auxiliary condition for a unique solution. However, the solution of the equation in strong form must be unique for certain time. The solution is regular in the beginning if it satisfies the above-mentioned compatibility condition. In this paper, the evolution of a smooth initial velocity field is calculated analytically using the differential form of the equation. Furthermore, a smooth example is given that violates the compatibility condition, therefore a regular solution does not exist at the initial time. 
For the sake of simplicity the domain is a rectangular cuboid, a channel configuration.
In that case, the steady-state solution is the well-known Poiseuille flow. Recently, \citet{Jozsa2019} obtained an analytical solution for controlled channel flow.

This paper is structured as follows. First, the analytic calculation method of the temporal evolution of the velocity is presented using the vorticity equation. A condition is obtained that can be used to investigate the regularity of the solution at the initial time. The benefit of this method is that the boundary conditions are necessary only for the velocity, since the pressure eliminated. This is followed by the discussion of the the compatibility condition. Then an analytic calculation example is given that violates the condition and the solution is irregular at the beginning. Finally, counterexamples are presented in which case the condition holds and the solution is regular.

%
\section{The calculation method of the analytic solution} \label{sec_method}
\subsection{Governing equations}
An incompressible Newtonian fluid can be described by the continuity equation
\begin{equation}\label{eq:continuity}
\pder{{u}_{i}}{x_i}=0
\end{equation}
and the Navier-Stokes equation in non-dimensional form
\begin{equation}\label{NS}
\pder{{u}_{i}}{t}+{u}_{j}\pder{{u}_{i}}{x_j}=-\pder{{p}}{x_i}+\frac{1}{\Rey}\pder{^2 {u}_i}{x^2_j},
\end{equation}
where ${u}_i$ is the non-dimensional velocity, ${p}$ is the non-dimensional pressure
, $\Rey$ is the Reynolds number, $i = \{1,2,3\}$ is a running index of the Einstein summation notation.
\begin{equation}\label{eq_Reynolds_number}
\Rey = \frac{U_0 h}{\nu},
\end{equation}
$U_0$ is the velocity scale, $h$ is the half gap of the channel, $\nu$ is the kinematic viscosity. The body forces are neglected.

The non-dimensional solution domain is a cuboid: $x_1\in[0, L_x], x_2\in[-1, 1], x_3\in[0, L_z]$ that represents a channel. %
%
The boundary conditions are the following. The velocity field is periodic in the $x_1$ and $x_3$ directions:
\begin{equation}\label{BC_periodic_1}
u_i(x_1, x_2, x_3,t)= u_i(x_1+L_x, x_2, x_3,t)
,
\end{equation}
\begin{equation}\label{BC_periodic_2}
u_i(x_1, x_2, x_3,t)= u_i(x_1, x_2, x_3+L_z,t).
\end{equation}
Usually, these directions are called streamwise and spanwise  directions, respectively. 
At $x_2=-1$ and $x_2=1$, the no-slip wall condition is prescribed, meaning that
\begin{equation}\label{BC_wall}
u_i(x_1, x_2=-1, x_3,t)=u_i(x_1, x_2=1, x_3,t)=0.
\end{equation}
$x_2$ is the wall-normal coordinate. 
The initial velocity field is given as
\begin{equation}\label{IC}
u_i(x_1,x_2,x_3, t=0) = u_{i,0}(x_1, x_2, x_3)
,
\end{equation}
which fulfills the continuity equation (\ref{eq:continuity}) and the boundary conditions. The initial velocity field example will be presented later.

\begin{figure}
	\begin{center}
		\includegraphics[width=6cm]{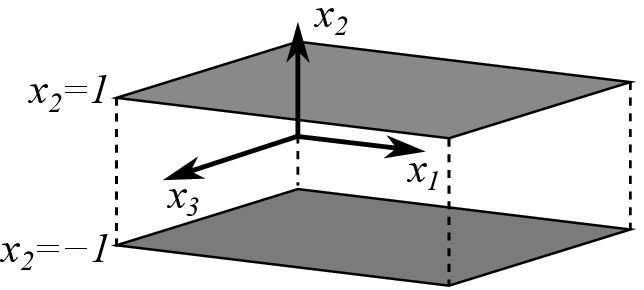} 
	\end{center}
	\caption{The channel geometry.}
	\label{fig_channel}
\end{figure}

%
%
A possible way to obtain the temporal derivative of the velocity field is using the curl of the Navier-Stokes equation, which known as vortex transport or vorticity equation, 
\begin{equation}\label{eq:vortex_transport}
\pder{{\omega}_{i}}{t}=-{u}_{j}\pder{{\omega}_{i}}{x_j}+{\omega}_{j}\pder{{u}_{i}}{x_j}+\frac{1}{\Rey}\pder{^2 {\omega}_i}{x_j^2},
\end{equation}
where $\omega_i$ is the vorticity, the curl of $u_i$. The benefit of this procedure is the elimination of the pressure.
$\omega_{i, 0}$ is the curl of the initial condition ($u_{i,0}$) that can be easily calculated for a given initial velocity. The temporal partial derivative of the vorticity at the initial time ($t=0$) can be determined using the vorticity equation (\ref{eq:vortex_transport})
\begin{equation}\label{eq:vortex_transport_initial}
\evalat{\pder{{\omega}_{i}}{t}}{t=0}=-{u}_{j,0}\pder{{\omega}_{i,0}}{x_j}+{\omega}_{j,0}\pder{{u}_{i,0}}{x_j}+\frac{1}{\Rey}\pder{^2 {\omega}_{i,0}}{x_j^2}
.
\end{equation}

Furthermore, the well-known vector algebra identity,
\begin{equation}\label{eq:vector_algebra_identitity}
\nabla\times\underbrace{(\nabla\times u_i)}_{\omega_i} = \nabla\underbrace{(\nabla \cdot u_i)}_{0}-\Delta u_i 
\end{equation}
and the symmetry of derivatives lead to 
\begin{equation}\label{eq:vector_algebra_identitity_2}
\nabla\times\left( \pder{\omega_i}{t}\right)  = -\Delta \left(\pder{u_i}{t}\right)
,
\end{equation}
since the temporal derivative of the velocity must be a solenoid field. $\nabla$ is the Nabla, $\Delta$ is the Laplace operator.
At the initial time, the temporal derivative of vorticity can be obtained from (\ref{eq:vortex_transport_initial}) and (\ref{eq:vector_algebra_identitity_2}) leads to three independent Poisson equations for the temporal derivative of each velocity component:
\begin{equation}\label{eq:Poisson_dudt}
\pder{^2}{x_j^2}\left(\evalat{\pder{{u}_{i}}{t}}{t=0}\right) = f_i(x_1, x_2,x_3)
,
\end{equation}
where
\begin{equation}\label{eq:ff}
f_k(x_1, x_2,x_3) = -\nabla\times\left( -{u}_{j,0}\pder{{\omega}_{i,0}}{x_j}+{\omega}_{j,0}\pder{{u}_{i,0}}{x_j}+\frac{1}{\Rey}\pder{^2 {\omega}_{i,0}}{x_j^2}\right) 
\end{equation}

The velocity boundary conditions of this problem are the same as those of the original problem, Eqs. (\ref{BC_periodic_1})-(\ref{BC_wall}). Since the Poisson equation is a linear differential equation, it has only one solution in the case of the Dirichlet boundary condition. However, it will be shown that for initial values that violate the compatibility condition, the solution of Eq. (\ref{eq:Poisson_dudt}) is not divergence-free. (It is mentioned that $f_k$ (\ref{eq:ff}) is undoubtedly divergence-free.) This unique solution of Eq. (\ref{eq:Poisson_dudt}) is therefore not acceptable, meaning that the solution of Eq. (\ref{NS}) must be irregular at $t=0$! 

\subsection{The compatibility condition}
\citet{Temam1982} proved (and later discussed \citep{Temam2006} the physical implications) that the solution is regular at the initial time, if and only if it fulfills the compatibility condition. In the case of the linear heat equation \citep{Temam2006}, the compatibility condition expresses that the initial value of the problem does not contradict the boundary condition. In contrast, the compatibility condition in the case of the Navier-Stokes equation is more complicated.
The condition is that the expression \citep{Temam1982}
\begin{equation}\label{eq_CC_Leray}
\frac{1}{\Rey}P(\Delta {u}_{i,0}) -  P\left({u}_{j,0}\pder{{u}_{i,0}}{x_j}\right)+g_i(0) = 0
\end{equation}
must be fulfilled at the boundaries; otherwise, the solution is irregular at $t=0$. (The sign of the second term is corrected here since it was wrongly positive in the cited paper.) $P$ is the Leray projector that gives the divergence-free part of a vector field. The projector can be expressed with the Helmholtz decomposition. $\Delta$ is the Laplace operator, $g_i$ is an arbitrary volumetric source term in the momentum equation, and is 0 in this study. 
\citet{Temam2006} proposed a calculation method to check the compatibility condition. First, the following equation is solved:  
\begin{equation}\label{eq_pressure}
\Delta {p}_0 = -\pder{{u}_{j,0}}{x_i}\pder{{u}_{i,0}}{x_j},
\end{equation}
which is the divergence of the momentum equation. It is also known as the Poisson equation for the pressure. (In the cited paper, the negative sign is omitted, which leads to the same condition, but in this form, $p$ is the widely used physical pressure.) The Neumann boundary condition for the pressure is
\begin{equation}\label{eq_pressure_BC}
\pder{{p}_{0}}{x_i}{n_i} = \left(\frac{1}{\Rey}\Delta{u}_{i,0}\right){n_i}
\end{equation}
in the case of a no-slip wall, which can be obtained from the momentum equation (\ref{NS}) at the stationary walls (\ref{BC_wall}) in the wall-normal direction. $n_i$ is the wall-normal vector. In the presented channel configuration, it is non-zero only for $i=2$.
\citet{Temam1982} states as a consequence of the compatibility condition after solving the problem (\ref{eq_pressure}) using (\ref{eq_pressure_BC}) that the temporal derivative of the velocity field is smooth at $t=0$, if the tangential components of $\pder{{p}_0}{x_i}$ are equal to the tangential components of $\frac{1}{\Rey}\Delta{u}_{i,0}$: 
\begin{equation}\label{eq_CC_cond}
\pder{{p}_{0}}{x_i}{t_i} = \left(\frac{1}{\Rey}\Delta{u}_{i,0}\right){t_i}
,
\end{equation}
where $t_i$ is an arbitrary vector perpendicular to $n_i$. The condition is actually the evaluation of the momentum equation at the wall in a tangential direction.
The physical meaning of Eq. (\ref{eq_CC_cond}) is that the pressure field prevents the tangential acceleration (the tangential movement) of the wall. The problem is that this additional condition cannot be prescribed for the Poisson equation (\ref{eq_pressure}), since then it would become overdetermined. This condition is a constraint for the initial velocity field. If it does not hold, there is no regular solution at the initial time which fulfills the stationary wall condition and the solenoid property at the same time. An irregular solution may exist but such a velocity field can exist only as a mathematical abstraction. It can smoothly evolve physically neither from a previous state nor to a next state. However, all this reasoning is valid only for a fully incompressible flow, which is an idealization. A real fluid is always compressible and the mathematical requirement that the pressure maintains, the incompressibility is relaxed.

%

\begin{figure}
	\centering
	\subfigure[]{
		\includegraphics[width=9cm]{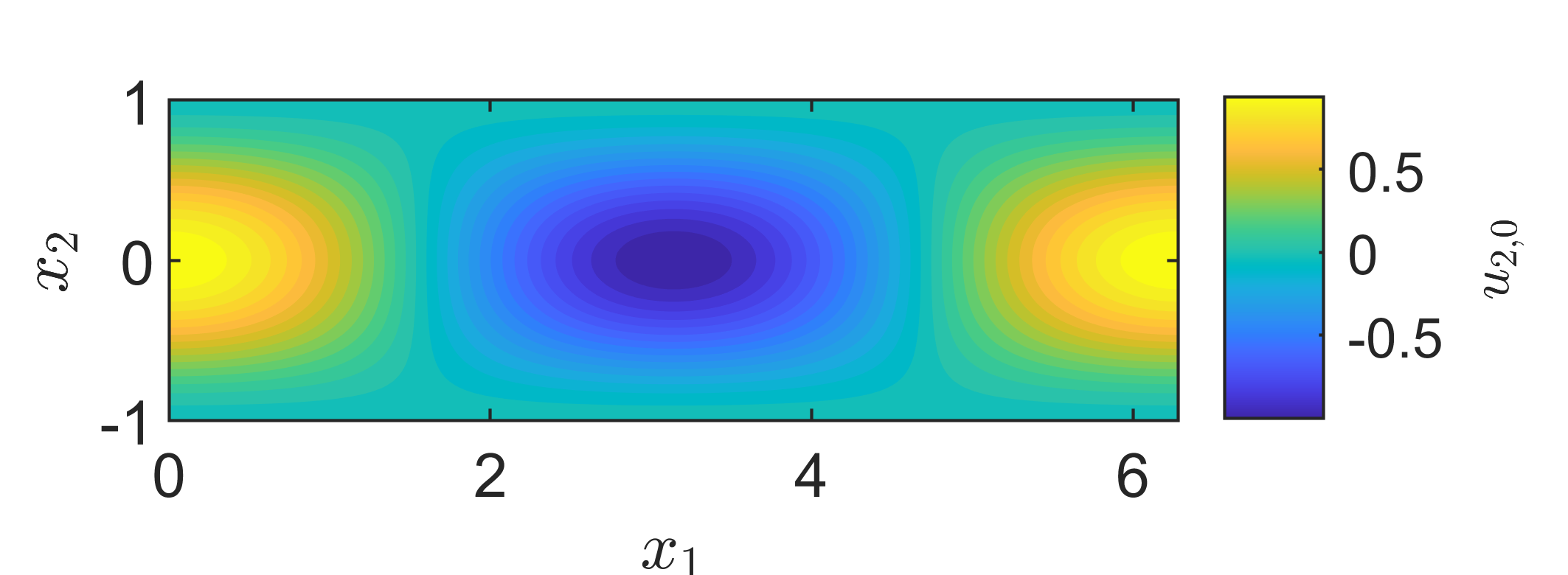}
		\label{fig_v}} 
	\subfigure[]{
		\includegraphics[width=9cm]{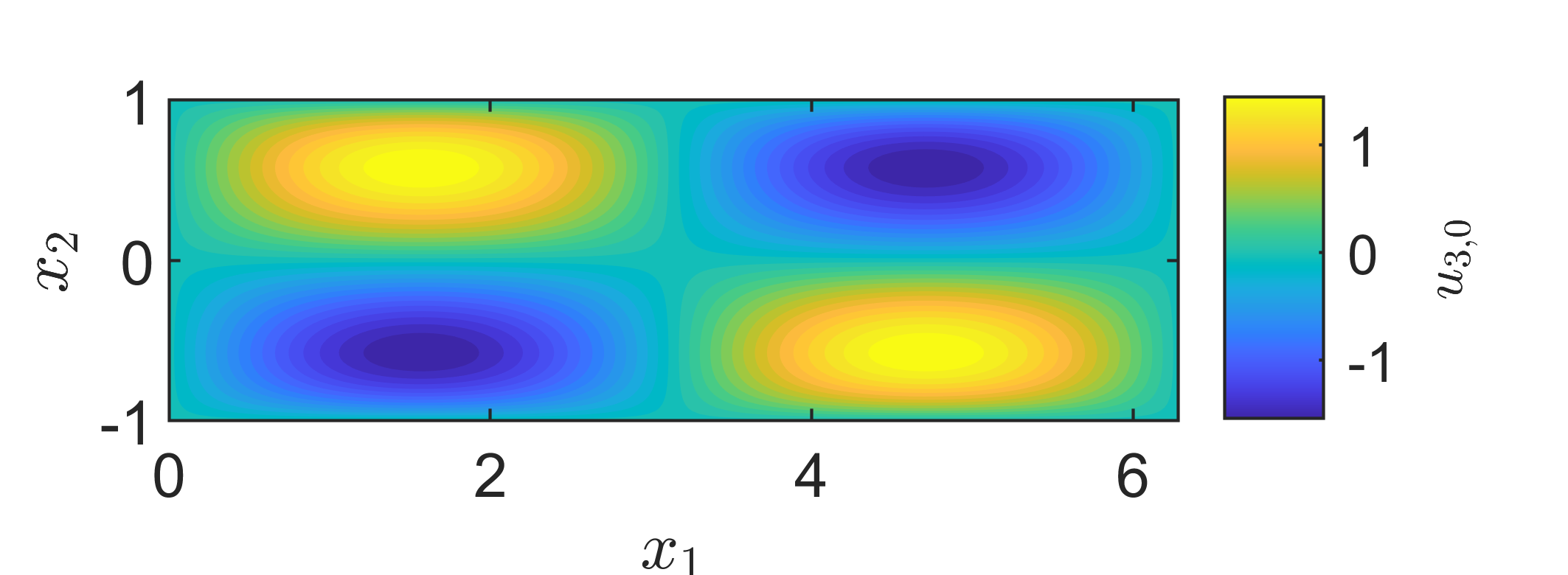}
		\label{fig_w}} 
	\caption{The slice ($z=0$) of the example velocity field. The parameters are $\alpha=1, \,\beta=1 ,\, Re = 80$ (a) $u_{2,0}$ (\ref{u_{i,0}_2b}); (b) $u_{3,0}$ Eq. (\ref{u_{i,0}_2c}).}
	\label{fig_vel}
\end{figure}

\begin{figure}
	\centering
	\includegraphics[width=9cm]{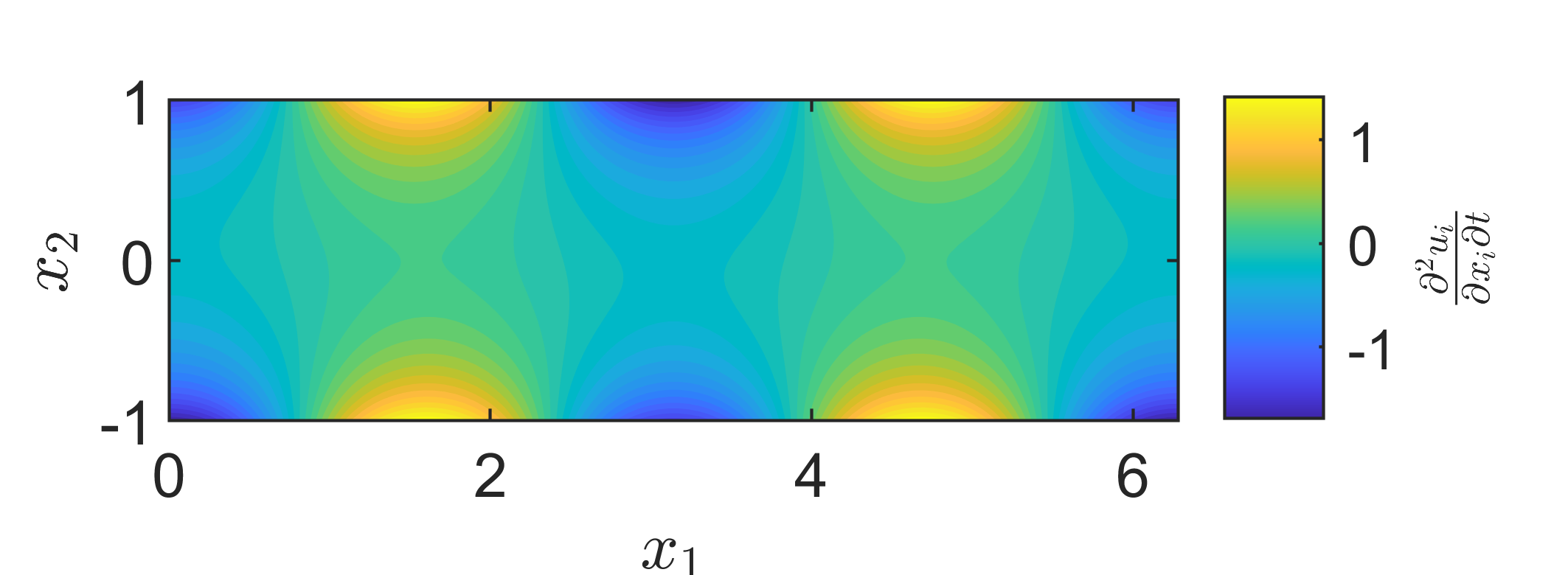}
	\caption{The divergence of $\evalat{\pder{{u}_{i}}{t}}{t=0}$ (Eq. \ref{eq:ut_i})  as the function of $x, y$  for the parameters $\alpha=1, \,\beta=1 ,\, Re = 80$ at $z=0$.}
	\label{fig:div_ut}
\end{figure}

Next, the temporal evolution of a given velocity field will be calculated analytically using Eq. (\ref{eq:Poisson_dudt}). The solution is assumed in a Fourier series, meaning that it is smooth in space. The only analytical solution will be shown to be wrong, meaning that a regular solution does not exist. It will be proven by evaluating the compatibility condition which does not hold.
In this paper, the calculation is presented for the following initial velocity field
\begin{align}\label{u_{i,0}_2}
u_{1,0}&=0\\
u_{2,0}&=\cos\left(\alpha \,x+\beta \,z\right)\,{\left(y^2-1\right)}^2\label{u_{i,0}_2b}\\
u_{3,0}&=-\frac{4\,y\,\sin\left(\alpha \,x+\beta \,z\right)\,\left(y^2-1\right)}{\beta }
\label{u_{i,0}_2c}
,
\end{align}
where $x=x_1,\; y=x_2,\; z=x_3,\; \alpha=2\pi/L_x,\; \beta=2\pi/L_z$. This is one of the simplest initial velocity fields, which is divergence-free and fulfills the boundary conditions. (Of course, $u_{1,0}$ is not necessarily 0. It was set to 0 for the sake of simplicity, since otherwise the analytical terms could become far too long for publication but this assumption does not restrict the message of the paper.) The velocity field is plotted in Figure \ref{fig_vel} at the plane $z=0$. The analytical calculation of the temporal derivative of the velocity, its divergence, and the compatibility condition are presented in the Appendix. 
The calculation is straightforward, but the long expressions would become even more lengthy in the case of a more complex initial field. Matlab 2019b Symbolic toolbox was used to reduce the solution time and the risk of miscalculation. 

The same calculation was repeated using numerical functions mainly based on the Chebyshev collocation method, similarly to the work of \citet{Falsaperla2019} and \citet{Nagy2022}. In this case, the parameters must have numerical values. The calculations were carried out multiple times for different parameters and the analytical and the numerical computation results were practically identical, meaning that no error was made in the calculation.

This initial velocity field does not fulfill the compatibility condition, and the corresponding problem has no regular solution at $t=0$. The divergence of the velocity derivative is shown in Figure \ref{fig:div_ut} for certain parameters, and its value is clearly non-zero. This means that the calculated temporal derivative of the velocity field is wrong, leading again to the conclusion that the solution must be irregular at the initial time. The calculated terms and their derivation are shifted to the Appendix because they are long expressions. Using the presented method, the calculated temporal derivative of the velocity field would lead to a solution that fulfills the boundary condition, but violates the divergence-free condition. It is possible to obtain a divergence-free solution using the pressure Poisson equation but it violates some component of the velocity boundary conditions at the wall. The two conditions cannot be fulfilled simultaneously that is the consequence of the the given initial velocity field violating the compatibility condition.  


Some further analysis was done to investigate, which similar initial fields fulfill the condition and lead therefore to a regular solution. The velocity field components have the form 
\begin{equation}\label{eq_polinomial_form}
u_{i,0} = a_{ui}(y)\cos\left(\alpha \,x+\beta \,z\right)+b_{ui}(y)\sin\left(\alpha \,x+\beta \,z\right)
\end{equation}
in all cases, where $a(y)$ and $b(y)$ are polynomial functions. In this form, the periodic boundary conditions (\ref{BC_periodic_1})-(\ref{BC_periodic_2}) are automatically fulfilled. The polynomials $a(y)$ and $b(y)$ must fulfill the wall boundary condition (\ref{BC_wall}) and the continuity equation (\ref{eq:continuity}). A possible choice for $a(y)$ and $b(y)$ is the form $c(y)=\tilde{c}(y)(y^2-1)$ for $u_{1,0}$ and $c(y)=\tilde{c}(y)(y^2-1)^2$ for $u_{2,0}$ to fulfill the boundary condition where $c(y)$ means $a(y)$ or $b(y)$ and $\tilde{c}(y)$ is an arbitrary polynomial. The last velocity component ($u_{3,0}$) cannot be arbitrary, it can be determined simply using the continuity equation (\ref{eq:continuity}).

After analyzing multiple initial velocity fields in the proposed form, the following conclusions can be drawn.
In the case of $u_{2,0}=0$, the solution of Equation (\ref{eq:Poisson_dudt}) is divergence-free in every case meaning that the solution is regular. 
In the case of $u_{2,0} \neq 0$, many attempts were made to obtain an analytical velocity field that fulfills the compatibility condition. The procedure failed since increasing the order of the polynomials leads to a very complicated analytical problem that cannot be solved by hand or on a personal computer. On the other hand, the equations can be easily solved numerically. In the following, an initial velocity field in the proposed polynomial and periodic form is searched which has a regular solution at $t=0$ and has wall-normal velocity components. 

First, the asymptotic eigenmodes of the corresponding Orr-Sommerfeld problem are used as initial conditions. They are the solution of the linearised Navier-Stokes equation where the base flow is the well-known Poiseuille flow. The eigenfunctions are obtained with the method of \citet{Juniper2013}, where the domain is discretized with $N=40$ Chebyshev polynomials. These eigenfunctions have the proposed form (\ref{eq_polinomial_form}) and usually have a wall-normal component. 
In this case, a spatially non-oscillating base flow, the Poiseuille flow is part of the initial velocity field, beside the wavelike velocity field in the proposed form (\ref{eq_polinomial_form}) that is usually called a perturbation. The smooth temporal evolution of the velocity field is known, if the amplitude of the perturbations is low. This means that the solution must be regular at $t=0$. It has been verified that these eigenvectors fulfill the compatibility condition and the divergence of the temporal derivative of the velocity field is zero. The question arises whether the compatibility condition holds for the non-linear Navier-Stokes equations too. After increasing the perturbation amplitude to reasonable finite values, the condition is numerically still fulfilled. These initial conditions therefore provide a regular solution for the Navier-Stokes problem.

Another approach was also followed to find initial velocity field that has wall-normal component and fulfills the compatibility condition. The proposed product of polynomial and trigonometric function form (\ref{eq_polinomial_form}) was used, the order of polynomial $\tilde{c}(y)$ was set to 4, meaning that $c(y)$ has the order of 6 or 8 depending on the velocity component. The same Chebyshev collocation method  was used to investigate the compatibility condition and the divergence of the temporal derivative of the velocity. The coefficients of the polynomials are varied to obtain a divergence-free temporal derivative of the velocity by using the non-linear optimization \textit{fsolve} function built in MATLAB. Multiple solutions were obtained; an example is given at $\alpha=1, \beta = 1, \Rey = 80$:
\begin{align}\label{eq_good_examplea}
a_{u1}&=\left(0.8147y^4+0.9058y^3+0.127y^2+0.9134y+0.6324\right)\left(y^2-1\right)\\
b_{u1}&=\left(0.09754y^4+0.2785y^3+0.5469y^2+0.9575y+0.9649\right)\left(y^2-1\right)\\
a_{u2}&=\left(-0.6011y^4-0.1986y^3+0.7068y^2+0.4689y+1.599\right)\left(y^2-1\right)^2\\
b_{u2}&=\left(0.1063y^4+0.2864y^3+0.8618y^2+0.3238y+1.537\right)\left(y^2-1\right)^2\label{eq_good_exampled}
\end{align}
the subscript of the Fourier coefficients ($a(y), b(y)$ polynomials) denotes the velocity component. The third velocity component can be easily obtained from the continuity equation. These functions (\ref{eq_good_examplea}-\ref{eq_good_exampled}) can be used to check numerical codes that cannot handle the problem of the irregularity of the solution at the initial time. In this way, other velocity fields can be obtained that have wall-normal component and still fulfill the compatibility condition.
Although these calculations are carried out using numerical techniques, there exist probably analytical initial velocity fields too, where $u_{2,0} \neq 0$, and the corresponding Navier-Stokes equation has a regular solution at $t=0$.


\section{Conclusion}\label{sec_conclusion}
The solution of the Navier-Stokes equation is a challenging problem. In the last century, mathematicians proved that a smooth solution exists for a finite time. However, in the case of physically relevant bounded domains, the solution is regular at the initial time if and only if the compatibility condition holds. In this paper, an analytic example is given in a channel flow, when the initial velocity field is smooth and divergence-free. At the same time, it does not fulfill the compatibility condition, and the corresponding Navier-Stokes equation has no regular solution at $t=0$. According to the authors' best knowledge, this is the first analytic example where the irregularity of the solution was demonstrated analytically. The problem is related to the presence of walls where the pressure Poisson equation can be overdetermined. Further single wave functions were investigated and all of them violated the conditions. Therefore, attempts are made to construct an initial velocity field which does not contradict the condition. In the absence of a wall-normal velocity component, the condition is always fulfilled in the presented configuration. Some further numerical calculations suggest that there are velocity fields with a non-zero wall-normal component that also fulfill the condition. An example is given in the paper. 
However, we have not been able to construct such a field analytically.

The results of the paper are direct consequences of Temam's\citep{Temam1982} mathematical proof about the regularity of the solution. Yet, the physical consequence is surprising. There exist many mathematically smooth velocity fields that fulfill the boundary conditions and the continuity equation, but they violate the governing equation in a non-trivial way. In these cases, the Navier-Stokes equation has no regular solution at the initial time and this may have a serious consequences for the initial part of time-dependent simulations.

The data that support the findings of this study are available from the corresponding author upon reasonable request.


\section*{Acknowledgments}
The work has been performed within the framework of the NKFI (National Research Development and Innovation Office of Hungary) project K124939 and K135436. 

\section*{Declaration of Competing Interest}
The authors declare that they have no known competing financial interests or personal relationships that could have appeared to influence the work reported in this paper.

\appendix
\section{The calculation of temporal derivative in an analytic example}

First, the curl of the initial velocity (\ref{u_{i,0}_2})-(\ref{u_{i,0}_2c}) is calculated:

\begin{align}\label{omega_{i,0}_2}
\omega_{1,0}&=\frac{\sin\left(\alpha \,x+\beta \,z\right)\,\left(\beta ^2\,y^4-2\,\beta ^2\,y^2+\beta ^2-12\,y^2+4\right)}{\beta }\\
\omega_{2,0}&=\frac{4\,\alpha \,y\,\cos\left(\alpha \,x+\beta \,z\right)\,\left(y^2-1\right)}{\beta }\\
\omega_{3,0}&=-\alpha \,\sin\left(\alpha \,x+\beta \,z\right)\,{\left(y^2-1\right)}^2
.
\end{align}

The temporal derivative of the vorticity ($\evalat{\pder{{\omega}_{i}}{t}}{t=0}$) is obtained from Eq. (\ref{eq:vortex_transport_initial}). These terms are not shown here.
Then the opposite and the curl of $\evalat{\pder{{\omega}_{i}}{t}}{t=0}$ is calculated according to Eq. (\ref{eq:ff}). It is convenient to express the term in a Fourier series in the following form: 
%
\begin{equation}\label{eq:f_i}
f_k = A_{k, 0}(y) +\sum_{j=1}^{2} \bigr\lbrace A_{k, j}(y) \cos\big(j(\alpha x+\beta z)\big) + B_{k, j}(y) \sin\big(j(\alpha x+\beta z)\big) \bigr\rbrace,
\end{equation}
where $A_{k, j}(y)$, $B_{k, j}(y)$ are polynomials. 

The coefficients are in the presented case:
\begin{align}\label{f_coeff_2}
A_{1,0}&=0\\
A_{1,1}&=0\\
B_{1,1}&=0\\
A_{1,2}&=0\\
B_{1,2}&=-8\,\alpha \,{\left(y^2-1\right)}^2\,\left(y^2+1\right)\\
A_{2,0}&=0\\
A_{2,1}&=-\frac{1}{\mathrm{Re}}\big(\left(-\alpha ^4-2\,\alpha ^2\,\beta ^2-\beta ^4\right)\,y^4+\left(2\,\alpha ^4+4\,\alpha ^2\,\beta ^2+24\,\alpha ^2+2\,\beta ^4+24\,\beta ^2\right)\,y^2\nonumber\\&-\alpha ^4-2\,\alpha ^2\,\beta ^2-8\,\alpha ^2-\beta ^4-8\,\beta ^2-24\big)\\
B_{2,1}&=0\\
A_{2,2}&=8\,y\,\left(3\,y^2+1\right)\,\left(y-1\right)\,\left(y+1\right)\\
B_{2,2}&=0\\
A_{3,0}&=0\\
A_{3,1}&=0\\
B_{3,1}&=\frac{4\,y\,\left(\alpha ^2+\beta ^2\right)\,\left(-\alpha ^2\,y^2+\alpha ^2-\beta ^2\,y^2+\beta ^2+12\right)}{\mathrm{Re}\,\beta }\\
A_{3,2}&=0\\
B_{3,2}&=\frac{4\,\left(2\,\alpha ^2\,y^6-2\,\alpha ^2\,y^4-2\,\alpha ^2\,y^2+2\,\alpha ^2-15\,y^4+6\,y^2+1\right)}{\beta }
\end{align}

The temporal derivative of the velocity must have a similar form:
\begin{equation}\label{eq:ut_i}
\evalat{\pder{{u}_{k}}{t}}{t=0} = a_{k, 0}(y) +\sum_{j=1}^{2} \bigr\lbrace a_{k, j}(y) \cos\big(j(\alpha x+\beta z)\big) + b_{k, j}(y) \sin\big(j(\alpha x+\beta z)\big) \bigr\rbrace.
\end{equation}
Its Laplacian is
\begin{align}\label{eq:ut_i_lap}
\pder{^2}{x_j^2}\left(\evalat{\pder{{u}_{k}}{t}}{t=0}\right) =& \tder{{}^2 a_{k, 0}(y)}{y^2} +\sum_{j=1}^{2} \left(-j^2(\alpha^2+\beta^2) a_{k, j}(y) + \tder{{}^2 a_{k, j}(y)}{y^2} \right) \cos\big(j(\alpha x+\beta z)\big) \nonumber\\&+\sum_{j=1}^{2} \left(-j^2(\alpha^2+\beta^2) b_{k, j}(y) + \tder{{}^2 b_{k, j}(y)}{y^2} \right) \sin\big(j(\alpha x+\beta z)\big) .
\end{align}
The solution of the Poisson Eq. (\ref{eq:Poisson_dudt}) for the temporal derivative of velocity is equivalent to solving (\ref{eq:f_i})=(\ref{eq:ut_i_lap}). The five terms (constant, sine, cosine, harmonics) in the three equations can be treated separately. The original Poisson equations become five independent 1D boundary value problems for each component of the vector field in the series expansion. For example: 
\begin{equation}
\tder{{}^2 a_{1, 0}(y)}{y^2}=A_{1, 0}(y)
\end{equation} 
or
\begin{equation}
\left(-\alpha^2 -\beta^2\right)a_{1, 1}(y)+\tder{{}^2 a_{1, 1}(y)}{y^2}=A_{1, 1}(y).
\end{equation} 

The periodic boundary conditions are automatically fulfilled because of the form of expressions. The no-slip boundary condition (\ref{BC_wall}) can be written as $a_{k,j}(y=1)=a_{k,j}(y=-1)=0$ for $j=\lbrace0,1,2\rbrace$ and similarly, $b_{k,j}(y=1)=b_{k,j}(y=-1)=0$ for $j=\lbrace1,2\rbrace$. The 15 ordinary differential equations are solved using Matlab symbolic and $a_{k, j}, b_{k, j}$ polynomials are obtained.
The solutions of the Poisson equations are:

\begin{align}\label{a_coeff_2}
a_{1,0}&=0\\
a_{1,1}&=0\\
b_{1,1}&=0\\
a_{1,2}&=0\\
b_{1,2}&=\frac{2\,\alpha \,y^6}{\alpha ^2+\beta ^2}\nonumber\\&-\frac{\alpha \,\left(-4\,\alpha ^6-12\,\alpha ^4\,\beta ^2+2\,\alpha ^4-12\,\alpha ^2\,\beta ^4+4\,\alpha ^2\,\beta ^2+6\,\alpha ^2-4\,\beta ^6+2\,\beta ^4+6\,\beta ^2-45\right)}{2\,{\left(\alpha ^2+\beta ^2\right)}^4}\nonumber\\&-\frac{\alpha \,y^4\,\left(2\,\alpha ^2+2\,\beta ^2-15\right)}{{\left(\alpha ^2+\beta ^2\right)}^2}-\frac{\alpha \,y^2\,\left(2\,\alpha ^4+4\,\alpha ^2\,\beta ^2+6\,\alpha ^2+2\,\beta ^4+6\,\beta ^2-45\right)}{{\left(\alpha ^2+\beta ^2\right)}^3}\nonumber\\&-\frac{\alpha \,{\mathrm{e}}^{-2\,y\,\sqrt{\alpha ^2+\beta ^2}}\,\left(16\,\alpha ^4+32\,\alpha ^2\,\beta ^2+84\,\alpha ^2+16\,\beta ^4+84\,\beta ^2+45\right)}{2\,\left({\mathrm{e}}^{-2\,\sqrt{\alpha ^2+\beta ^2}}+{\mathrm{e}}^{2\,\sqrt{\alpha ^2+\beta ^2}}\right)\,{\left(\alpha ^2+\beta ^2\right)}^4}\nonumber\\&-\frac{\alpha \,{\mathrm{e}}^{2\,y\,\sqrt{\alpha ^2+\beta ^2}}\,\left(16\,\alpha ^4+32\,\alpha ^2\,\beta ^2+84\,\alpha ^2+16\,\beta ^4+84\,\beta ^2+45\right)}{2\,\left({\mathrm{e}}^{-2\,\sqrt{\alpha ^2+\beta ^2}}+{\mathrm{e}}^{2\,\sqrt{\alpha ^2+\beta ^2}}\right)\,{\left(\alpha ^2+\beta ^2\right)}^4}
\end{align}
\begin{align}
a_{2,0}&=0\\
a_{2,1}&=\frac{2\,y^2\,\left(\alpha ^2+\beta ^2+6\right)}{\mathrm{Re}}-\frac{8\,{\mathrm{e}}^{\left(y+1\right)\,\sqrt{\alpha ^2+\beta ^2}}}{\mathrm{Re}\,\left({\mathrm{e}}^{2\,\sqrt{\alpha ^2+\beta ^2}}+1\right)}-\frac{y^4\,\left(\alpha ^2+\beta ^2\right)}{\mathrm{Re}}-\frac{\alpha ^2+\beta ^2+4}{\mathrm{Re}}\nonumber\\&-\frac{8\,{\mathrm{e}}^{-y\,\sqrt{\alpha ^2+\beta ^2}}}{\mathrm{Re}\,\left({\mathrm{e}}^{-\sqrt{\alpha ^2+\beta ^2}}+{\mathrm{e}}^{\sqrt{\alpha ^2+\beta ^2}}\right)}\\
b_{2,1}&=0\\
a_{2,2}&=\frac{2\,y^3\,\left(2\,\alpha ^2+2\,\beta ^2-15\right)}{{\left(\alpha ^2+\beta ^2\right)}^2}-\frac{6\,y^5}{\alpha ^2+\beta ^2}+\frac{y\,\left(2\,\alpha ^4+4\,\alpha ^2\,\beta ^2+6\,\alpha ^2+2\,\beta ^4+6\,\beta ^2-45\right)}{{\left(\alpha ^2+\beta ^2\right)}^3}\nonumber\\&+\frac{3\,{\mathrm{e}}^{-2\,y\,\sqrt{\alpha ^2+\beta ^2}}\,\left(8\,\alpha ^2+8\,\beta ^2+15\right)}{\left({\mathrm{e}}^{-2\,\sqrt{\alpha ^2+\beta ^2}}-{\mathrm{e}}^{2\,\sqrt{\alpha ^2+\beta ^2}}\right)\,{\left(\alpha ^2+\beta ^2\right)}^3}-\frac{3\,{\mathrm{e}}^{2\,y\,\sqrt{\alpha ^2+\beta ^2}}\,\left(8\,\alpha ^2+8\,\beta ^2+15\right)}{\left({\mathrm{e}}^{-2\,\sqrt{\alpha ^2+\beta ^2}}-{\mathrm{e}}^{2\,\sqrt{\alpha ^2+\beta ^2}}\right)\,{\left(\alpha ^2+\beta ^2\right)}^3}\\
b_{2,2}&=0
\end{align}

\begin{align}
a_{3,0}&=0\\
a_{3,1}&=\frac{2\,y^2\,\left(\alpha ^2+\beta ^2+6\right)}{\mathrm{Re}}-\frac{8\,{\mathrm{e}}^{\left(y+1\right)\,\sqrt{\alpha ^2+\beta ^2}}}{\mathrm{Re}\,\left({\mathrm{e}}^{2\,\sqrt{\alpha ^2+\beta ^2}}+1\right)}-\frac{y^4\,\left(\alpha ^2+\beta ^2\right)}{\mathrm{Re}}-\frac{\alpha ^2+\beta ^2+4}{\mathrm{Re}}\nonumber\\&-\frac{8\,{\mathrm{e}}^{-y\,\sqrt{\alpha ^2+\beta ^2}}}{\mathrm{Re}\,\left({\mathrm{e}}^{-\sqrt{\alpha ^2+\beta ^2}}+{\mathrm{e}}^{\sqrt{\alpha ^2+\beta ^2}}\right)}\\
b_{3,1}&=0\\
a_{3,2}&=\frac{2\,y^3\,\left(2\,\alpha ^2+2\,\beta ^2-15\right)}{{\left(\alpha ^2+\beta ^2\right)}^2}-\frac{6\,y^5}{\alpha ^2+\beta ^2}+\frac{y\,\left(2\,\alpha ^4+4\,\alpha ^2\,\beta ^2+6\,\alpha ^2+2\,\beta ^4+6\,\beta ^2-45\right)}{{\left(\alpha ^2+\beta ^2\right)}^3}\nonumber\\&+\frac{3\,{\mathrm{e}}^{-2\,y\,\sqrt{\alpha ^2+\beta ^2}}\,\left(8\,\alpha ^2+8\,\beta ^2+15\right)}{\left({\mathrm{e}}^{-2\,\sqrt{\alpha ^2+\beta ^2}}-{\mathrm{e}}^{2\,\sqrt{\alpha ^2+\beta ^2}}\right)\,{\left(\alpha ^2+\beta ^2\right)}^3}-\frac{3\,{\mathrm{e}}^{2\,y\,\sqrt{\alpha ^2+\beta ^2}}\,\left(8\,\alpha ^2+8\,\beta ^2+15\right)}{\left({\mathrm{e}}^{-2\,\sqrt{\alpha ^2+\beta ^2}}-{\mathrm{e}}^{2\,\sqrt{\alpha ^2+\beta ^2}}\right)\,{\left(\alpha ^2+\beta ^2\right)}^3}\\
b_{3,2}&=0.
\end{align}

The divergence of $\evalat{\pder{{u}_{k}}{t}}{t=0}$ (Eq. \ref{eq:ut_i}) is taken, its Fourier coefficents are:
\begin{align}\label{d_coeff_2}
da_{0}&=0\\
da_{1}&=\frac{8\,{\mathrm{e}}^{-y\,\sqrt{\alpha ^2+\beta ^2}}\,\left({\mathrm{e}}^{\left(2\,y+1\right)\,\sqrt{\alpha ^2+\beta ^2}}-{\mathrm{e}}^{\sqrt{\alpha ^2+\beta ^2}}\right)\,C_1}{\mathrm{Re}\,\left({\mathrm{e}}^{4\,\sqrt{\alpha ^2+\beta ^2}}-1\right)}\\
C_1 &=3\,{\mathrm{e}}^{2\,\sqrt{\alpha ^2+\beta ^2}}+\sqrt{\alpha ^2+\beta ^2}-{\mathrm{e}}^{2\,\sqrt{\alpha ^2+\beta ^2}}\,\sqrt{\alpha ^2+\beta ^2}+3\\
db_{1}&=0
\end{align}

\begin{align}
da_{2}&=-\frac{{\mathrm{e}}^{-2\,\left(y-1\right)\,\sqrt{\alpha ^2+\beta ^2}}\,\left({\mathrm{e}}^{4\,y\,\sqrt{\alpha ^2+\beta ^2}}+1\right) }{\left({\mathrm{e}}^{8\,\sqrt{\alpha ^2+\beta ^2}}-1\right)\,{\left(\alpha ^2+\beta ^2\right)}^4} C_2\\
C_2 &= \left(16\,{\mathrm{e}}^{4\,\sqrt{\alpha ^2+\beta ^2}}-16\right)\,{\left(\alpha ^2+\beta ^2\right)}^3-\left(90\,{\mathrm{e}}^{4\,\sqrt{\alpha ^2+\beta ^2}}+90\right)\,{\left(\alpha ^2+\beta ^2\right)}^{3/2}\nonumber\\
&+\alpha ^2\,\left(45\,{\mathrm{e}}^{4\,\sqrt{\alpha ^2+\beta ^2}}-45\right)+\alpha ^4\,\left(84\,{\mathrm{e}}^{4\,\sqrt{\alpha ^2+\beta ^2}}-84\right)+\beta ^2\,\left(45\,{\mathrm{e}}^{4\,\sqrt{\alpha ^2+\beta ^2}}-45\right)\nonumber\\
&+\beta ^4\,\left(84\,{\mathrm{e}}^{4\,\sqrt{\alpha ^2+\beta ^2}}-84\right)-\alpha ^2\,\left(48\,{\mathrm{e}}^{4\,\sqrt{\alpha ^2+\beta ^2}}+48\right)\,{\left(\alpha ^2+\beta ^2\right)}^{3/2}\nonumber\\
&-\beta ^2\,\left(48\,{\mathrm{e}}^{4\,\sqrt{\alpha ^2+\beta ^2}}+48\right)\,{\left(\alpha ^2+\beta ^2\right)}^{3/2}+\alpha ^2\,\beta ^2\,\left(168\,{\mathrm{e}}^{4\,\sqrt{\alpha ^2+\beta ^2}}-168\right)\\
db_{2}&=0
\end{align}

For the parameters $\alpha=1, \,\beta=1 ,\, Re = 80$  the divergence was plotted as the function of $x, y$ at $z=0$ in Fig. \ref{fig:div_ut}. 

The compatibility condition is obtained similarly. The pressure Poisson eq. (\ref{eq_pressure}) is solved with the boundary condition (\ref{eq_pressure_BC}) for $i=2$ for each periodic component.The compatibility condition (\ref{eq_CC_cond}) is evaluated at the walls ($y=\pm1$) in only the two main tangential directions: the streamwise and spanwise directions. Here, the difference between the two sides of Eq. (\ref{eq_CC_cond}) is presented with the Fourier coefficients. Due to the symmetry property of the example, the Fourier coefficients at the top and bottom walls are almost the same; only the signs are different in some components.  

The Fourier coefficients of the difference at the top ($y=1$) evaluating the compatibility condition in the streamwise ($x$) direction:
\begin{align}\label{diffx_coeff_2}
CCxa_{0}&=0\\
CCxa_{1}&=0\\
CCxb_{1}&=\frac{8\,\alpha \,\left({\mathrm{e}}^{2\,\sqrt{\alpha ^2+\beta ^2}}-1\right)}{\mathrm{Re}\,\left({\mathrm{e}}^{2\,\sqrt{\alpha ^2+\beta ^2}}+1\right)\,\sqrt{\alpha ^2+\beta ^2}}\\
CCxa_{2}&=0\\
CCxb_{2}&=-2\,\alpha \, C_3\\
C_3&=\frac{2\,\alpha ^2+2\,\beta ^2-15}{2\,{\left(\alpha ^2+\beta ^2\right)}^2}-\frac{1}{\alpha ^2+\beta ^2}+\frac{2\,\alpha ^4+4\,\alpha ^2\,\beta ^2+6\,\alpha ^2+2\,\beta ^4+6\,\beta ^2-45}{2\,{\left(\alpha ^2+\beta ^2\right)}^3}+
\\&+\frac{-4\,\alpha ^6-12\,\alpha ^4\,\beta ^2+2\,\alpha ^4-12\,\alpha ^2\,\beta ^4+4\,\alpha ^2\,\beta ^2+6\,\alpha ^2-4\,\beta ^6+2\,\beta ^4+6\,\beta ^2-45}{4\,{\left(\alpha ^2+\beta ^2\right)}^4}+
\\&+\frac{24\,\alpha ^2+24\,\beta ^2+45}{2\,\left({\mathrm{e}}^{4\,\sqrt{\alpha ^2+\beta ^2}}-1\right)\,{\left(\alpha ^2+\beta ^2\right)}^{7/2}}+\frac{3\,{\mathrm{e}}^{4\,\sqrt{\alpha ^2+\beta ^2}}\,\left(8\,\alpha ^2+8\,\beta ^2+15\right)}{2\,\left({\mathrm{e}}^{4\,\sqrt{\alpha ^2+\beta ^2}}-1\right)\,{\left(\alpha ^2+\beta ^2\right)}^{7/2}}
\end{align}

The difference coefficients at the top ($y=1$) evaluating the compatibility condition in the spanwise ($z$) direction:
\begin{align}\label{diffz_coeff_2}
CCza_{0}&=0\\
CCza_{1}&=0\\
CCzb_{1}&=\frac{8\,\beta \,\left({\mathrm{e}}^{2\,\sqrt{\alpha ^2+\beta ^2}}-1\right)}{\mathrm{Re}\,\left({\mathrm{e}}^{2\,\sqrt{\alpha ^2+\beta ^2}}+1\right)\,\sqrt{\alpha ^2+\beta ^2}}-\frac{24}{\mathrm{Re}\,\beta }\\
CCza_{2}&=0\\
CCzb_{2}&=-2\,\beta \, C_3
\end{align}
The coefficients are zero if the compatibility condition is fulfilled. However, it is not the case here.

\bibliographystyle{unsrtnat}
\bibliography{biblography_V1}

\begin{thebibliography}{13}
\providecommand{\natexlab}[1]{#1}
\providecommand{\url}[1]{\texttt{#1}}
\expandafter\ifx\csname urlstyle\endcsname\relax
  \providecommand{\doi}[1]{doi: #1}\else
  \providecommand{\doi}{doi: \begingroup \urlstyle{rm}\Url}\fi

\bibitem[Ladyzhenskaia(1969)]{Ladyzhenskaia1969}
O.~A. Ladyzhenskaia.
\newblock \emph{The Mathematical Theory of Viscous Incompressible Flow, 2nd
  ed.}
\newblock Gordon and Breach, New York, 1969.

\bibitem[Temam(1982)]{Temam1982}
R.~Temam.
\newblock Behaviour at time t = 0 of the solutions of semi-linear evolution
  equations.
\newblock \emph{Journal of Differential Equations}, 43\penalty0 (1):\penalty0
  73--92, 1982.
\newblock ISSN 0022-0396.
\newblock \doi{https://doi.org/10.1016/0022-0396(82)90075-4}.
\newblock URL
  \url{https://www.sciencedirect.com/science/article/pii/0022039682900754}.

\bibitem[Temam(2006)]{Temam2006}
Roger Temam.
\newblock Suitable initial conditions.
\newblock \emph{J. Comput. Phys.}, 218\penalty0 (2):\penalty0 443–450,
  November 2006.
\newblock ISSN 0021-9991.
\newblock \doi{10.1016/j.jcp.2006.03.033}.
\newblock URL \url{https://doi.org/10.1016/j.jcp.2006.03.033}.

\bibitem[Orszag and Israeli(1974)]{Orszag1974}
S~A Orszag and M~Israeli.
\newblock Numerical simulation of viscous incompressible flows.
\newblock \emph{Annual Review of Fluid Mechanics}, 6\penalty0 (1):\penalty0
  281--318, 1974.
\newblock \doi{10.1146/annurev.fl.06.010174.001433}.

\bibitem[Gresho and Sani(1987)]{Gresho1987}
Philip~M. Gresho and Robert~L. Sani.
\newblock On pressure boundary conditions for the incompressible navier-stokes
  equations.
\newblock \emph{International Journal for Numerical Methods in Fluids},
  7\penalty0 (10):\penalty0 1111--1145, 1987.
\newblock \doi{https://doi.org/10.1002/fld.1650071008}.
\newblock URL
  \url{https://onlinelibrary.wiley.com/doi/abs/10.1002/fld.1650071008}.

\bibitem[Johnston and Liu(2004)]{Johnston2004}
Hans Johnston and Jian-Guo Liu.
\newblock Accurate, stable and efficient navier–stokes solvers based on
  explicit treatment of the pressure term.
\newblock \emph{Journal of Computational Physics}, 199\penalty0 (1):\penalty0
  221--259, 2004.
\newblock ISSN 0021-9991.
\newblock \doi{https://doi.org/10.1016/j.jcp.2004.02.009}.
\newblock URL
  \url{https://www.sciencedirect.com/science/article/pii/S002199910400083X}.

\bibitem[Gallavotti(2002)]{Gallavotti2002}
Giovanni Gallavotti.
\newblock \emph{Foundations of Fluid Dynamics}.
\newblock Springer-Verlag Berlin Heidelberg, 2002.
\newblock \doi{10.1007/978-3-662-04670-8}.

\bibitem[Orr(1907)]{Orr1907}
William~McF Orr.
\newblock {The Stability or Instability of the Steady Motions of a Perfect
  Liquid and of a Viscous Liquid. Part II: A Viscous Liquid}.
\newblock \emph{Proc. R. Irish Acad.}, 27:\penalty0 69--138, 1907.
\newblock ISSN 19454589.
\newblock \doi{10.1017/CBO9781107415324.004}.

\bibitem[Heywood(1980)]{Heywood1980}
John~G. Heywood.
\newblock Auxiliary flux and pressure conditions for navier-stokes problems.
\newblock In Reimund Rautmann, editor, \emph{Approximation Methods for
  Navier-Stokes Problems}, pages 223--234, Berlin, Heidelberg, 1980. Springer
  Berlin Heidelberg.
\newblock ISBN 978-3-540-38550-9.

\bibitem[Józsa(2019)]{Jozsa2019}
Tamás~István Józsa.
\newblock Analytical solutions of incompressible laminar channel and pipe flows
  driven by in-plane wall oscillations.
\newblock \emph{Physics of Fluids}, 31\penalty0 (8):\penalty0 083605, 2019.
\newblock \doi{10.1063/1.5104356}.
\newblock URL \url{https://doi.org/10.1063/1.5104356}.

\bibitem[Falsaperla et~al.(2019)Falsaperla, Giacobbe, and
  Mulone]{Falsaperla2019}
Paolo Falsaperla, Andrea Giacobbe, and Giuseppe Mulone.
\newblock Nonlinear stability results for plane couette and poiseuille flows.
\newblock \emph{Phys. Rev. E}, 100:\penalty0 013113, Jul 2019.
\newblock \doi{10.1103/PhysRevE.100.013113}.
\newblock URL \url{https://link.aps.org/doi/10.1103/PhysRevE.100.013113}.

\bibitem[Nagy(2022)]{Nagy2022}
Péter~Tamás Nagy.
\newblock Enstrophy change of the reynolds-orr solution in channel flow.
\newblock \emph{Phys. Rev. E}, 105:\penalty0 035108, 2022.
\newblock \doi{https://doi.org/10.1103/PhysRevE.105.035108}.

\bibitem[Juniper et~al.(2013)Juniper, Hanifi, and Theofilis]{Juniper2013}
Matthew~P. Juniper, Ardeshir Hanifi, and Vassilios Theofilis.
\newblock Modal stability theory: Lecture notes from the flow-nordita summer
  school on advanced instability methods for complex flows, stockholm, sweden,
  2013.
\newblock \emph{Appl. Mech. Rev}, 66\penalty0 (2), September 2013.
\newblock ISSN 0003-6900.
\newblock URL \url{https://doi.org/10.1115/1.4026604}.

\end{thebibliography}



\end{document}